\documentclass[11pt]{amsart}
\usepackage{amssymb,amsmath}
\usepackage[T1]{fontenc}

\theoremstyle{plain}
\newtheorem{thm}{Th\'eor\`eme}[section]
\newtheorem{cor}[thm]{Corollaire}
\newtheorem{pro}[thm]{Proposition}
\newtheorem{lem}[thm]{Lemme}

\theoremstyle{definition}
\newtheorem*{defi}{D\'efinition}
\newtheorem*{eg}{Exemple}
\newtheorem*{egs}{Exemples}
\newtheorem{rem}[thm]{Remarque}

\def\vs{\vspace{0.3cm}}
\def\hs{\hspace{0.1cm}}

\def\C{\mathbb{C}}
\def\P{\mathbb{P}}
\def\Z{\mathbb{Z}}
\def\R{\mathbb{R}}
\def\Bir{\mathsf{Bir}}
\def\Aut{\mathsf{Aut}}
\def\PGL{\mathsf{PGL}}
\def\deg{\mathsf{deg}\,}

\def\og{\leavevmode\raise.3ex\hbox{$\scriptscriptstyle\langle\!\langle$~}}
\def\fg{\leavevmode\raise.3ex\hbox{~$\!\scriptscriptstyle\,\rangle\!\rangle$}}

\addtocounter{section}{0}             % Start with section 1
\numberwithin{equation}{section}       % Number formulas within sections
\begin{document}

\title{Sur les sous-groupes nilpotents du groupe de \textsc{Cremona}}
\author{Julie \textsc{D\'eserti}}

\email{julie.deserti@univ-rennes1.fr}
\address{IRMAR, UMR 6625 du CNRS\\
         Universit\'e de Rennes I\\
         35042 Rennes, France}

\begin{abstract}
We describe the nilpotent subgroups of the group $\Bir(\P^2
(\C))$ of birational transformations of the complex projective plane.
Let $\mathsf{N}$ be a nilpotent subgroup of class $k>1$;
then either each element of $\mathsf{N}$ has finite order,
or $\mathsf{N}$ is virtually metabelian. 
\end{abstract}

\maketitle

\section{Introduction}

Dans \cite{[Gh]} \textsc{Ghys} montre que {\it tout
sous-groupe nilpotent de
$\mathsf{Diff}^\omega(\mathbb{S}^2)$ est m\'etab\'elien} et en
d\'eduit que {\it si $\Gamma$ est un sous-groupe d'indice fini de
$\mathsf{SL}_n(\Z)$, avec $n\geq 4$, alors tout morphisme de
$\Gamma$ dans $\mathsf{Diff}^\omega(\mathbb{S}^2)$ est d'image
finie}. Nous \'etudions ici, via des techniques compl\`etement
diff\'erentes, les sous-groupes nilpotents du groupe des
transformations birationnelles du plan projectif
complexe~; ce groupe, que nous noterons $\Bir(
\P^2(\C))$, est aussi appel\'e groupe de \textsc{Cremona}.

Nous dirons qu'un groupe est \textsf{nilpotent fortement de
longueur $k$} (n. f. de longueur~$k$) s'il est nilpotent de
longueur $k$ et non virtuellement de longueur $k-1$. Rappelons que
le groupe de \textsc{Jonqui\`eres} est le sous-groupe des
transformations birationnelles qui pr\'eservent la fibration
rationnelle $y=$ cte (\emph{voir} \ref{jon}).

\begin{thm}\label{nilpotent}
Soient $\mathsf{N}$ un groupe n. f. de longueur $k>1$
et $\rho$ un morphisme injectif de $\mathsf{N}$ dans
$\Bir(\P^2(\C))$. Le groupe $\mathsf{N}$ v\'erifie l'une
des propri\'et\'es suivantes~:
\begin{itemize}
\item $\mathsf{N}$ est virtuellement m\'etab\'elien~;

\item $\mathsf{N}$ est de torsion.
\end{itemize}
\end{thm}

\begin{cor}
Soit $\mathsf{G}$ un groupe contenant un sous-groupe 
n. f. de longueur $k>1$ sans torsion et non 
virtuellement m\'etab\'elien. Alors il n'existe pas de
repr\'esentation fid\`ele de $\mathsf{G}$ dans le
groupe de \textsc{Cremona}.
\end{cor}

\noindent En particulier nous avons, dans l'esprit de 
\cite{[Zi]}, la cons\'equence suivante~:
{\it soit $\Gamma$ un sous-groupe d'indice fini de $\mathsf{SL}_n(
\Z)$~; d\`es que $n\geq 5$ le groupe $\Gamma$ ne se plonge
pas dans $\Bir(\P^2(\C)).$} Notons que nous avions obtenu 
un r\'esultat plus pr\'ecis ($n\geq 4$) dans \cite{[De]}
mais en utilisant des techniques plus complexes.

\section{Quelques rappels}

\subsection{Sur le groupe de \textsc{Jonqui\`eres}}\label{jon}
Si $f_S$ est une famille de transformations birationnelles,
alors $\langle f_S\rangle$ est le groupe engendr\'e par la
famille $f_S$. Dans
une carte affine $(x,y)$ de $\P^2(\C)$, si $f$ est un \'el\'ement
de $\Bir(\P^2(\C))$ nous notons $f$ par ses deux composantes
$(f_1(x,y),f_2(x,y))$.

Tout groupe de transformations birationnelles qui pr\'eserve
une fibration rationnelle est isomorphe au \textsf{groupe 
$\mathsf{J}$ de \textsc{Jonqui\`eres}}, {\it i.e.} isomorphe
\`a $\PGL_2(\C(y))\rtimes\PGL_2(\C)$.

\subsection{Dynamique et distorsion}

Soit $X$ une surface complexe compacte. Le {\sf premier degr\'e
dynamique} d'une transformation birationnelle $f\hs \colon\hs
X\dasharrow X$ est d\'efini par~: $$\lambda(f):= \limsup_{n\to
+\infty} |(f^n)^*|^{1/n}$$ o\`u $f^*$ d\'esigne l'application
lin\'eaire induite par $f$ de $ \mathsf{H}^{1,1}(X,\R)$ dans
lui-m\^eme (\cite{[Di-Fa]}) et $|\hs .\hs|$ une norme
sur $\mathsf{End}( \mathsf{H}^{1,1}(X,\R)).$ Ce nombre est
minor\'e par $1$ (\emph{voir} \cite{[Fr], [RS]}).

\begin{defi}\label{virtid}
Soit $S$ une surface complexe compacte. La transformation
birationnelle $f\hs\colon\hs S\dasharrow S$ est dite {\sf
virtuellement isotope \`a l'identit\'e} s'il existe une
surface $X$, une transformation birationnelle $\eta\hs
\colon\hs S\dasharrow X$ et un entier $k>0$ tels que $\eta
f^k\eta^{-1}$ soit un automorphisme de $X$ isotope \`a l'identit\'e.

\noindent Deux transformations birationnelles $f$ et $g$
sur $S$ sont dites {\sf simultan\'ement virtuellement isotopes
\`a l'identit\'e} si le couple $(\eta,X)$ est commun
\`a $f$ et $g.$
\end{defi}

Rappelons le r\'esultat suivant d\^u \`a \textsc{Diller} et
\textsc{Favre}.

\begin{thm}[\cite{[Di-Fa]}, th\'eor\`eme~0.2]\label{dillerfavre}
Soit $f$ une transformation birationnelle. Supposons que
$\lambda(f)=1~;$ alors $f$ satisfait une et une seule des
conditions sui\-vantes~:

\begin{itemize}
\item la suite $|(f^n)^*|$ est born\'ee et $f$ est
virtuellement isotope \`a l'identit\'e~;

\item la suite $|(f^n)^*|$ est \`a croissance lin\'eaire,
$f$ laisse une unique fibration rationnelle invariante et
n'est pas conjugu\'ee \`a un automorphisme~;

\item la suite $|(f^n)^*|$ est \`a croissance quadratique et
$f$ est, \`a conjugaison birationnelle pr\`es, un automorphisme
qui pr\'eserve une unique fibration elliptique.
\end{itemize}
\end{thm}

\begin{defi}
Soient $\mathsf{G}$ un groupe de type fini, $\{a_1$, $\ldots$, $a_n\}$ une
partie g\'en\'eratrice de $\mathsf{G}$ et $f$ un \'el\'ement de $\mathsf{G}.$

\begin{itemize}

\item[{\sl (i)}] La {\sf longueur} de $f$, not\'ee $\|f\|$, est le
plus petit entier $k$ pour lequel il existe une suite $(s_1,
\ldots,s_k)$ d'\'el\'ements de $\{a_1,\ldots,a_n,a_1^{-1},
\ldots,a_n^{-1}\}$ telle que $f=s_1\ldots s_k$. \vs

\item[{\sl (ii)}] La quantit\'e $\displaystyle\lim_{k \to +
\infty}\|f^k\|/k$ est la {\sf longueur stable} de $f$ (\emph{voir}
\cite{[Ha]}). \vs

\item[{\sl (iii)}] Un \'el\'ement $f$ de $\mathsf{G}$ est {\sf distordu}
s'il est d'ordre infini et si sa longueur stable est
nulle.
\end{itemize}
\end{defi}

\begin{rem}\label{dist}
Soit $\langle \mathsf{f},\mathsf{g},\mathsf{h}\hs|\hs[\mathsf{f},
\mathsf{h}]=[\mathsf{g},\mathsf{h}]=\mathsf{id},\hs[\mathsf{f},
\mathsf{g}]=\mathsf{h}\rangle$ un groupe de \textsc{Heisenberg}~;
le g\'en\'erateur $\mathsf{h}$ est distordu.
\end{rem}

\begin{lem}[\cite{[De]}]\label{degdyn}
Soient $f$ un \'el\'ement distordu d'un groupe de type fini
$\mathsf{G}$ et $\tau$ un morphisme de $\mathsf{G}$ dans
$\Bir(\P^2(\C))$. Le premier degr\'e dynamique de $\tau(f)$
vaut $1$.
\end{lem}

\subsection{Automorphismes isotopes \`a l'identit\'e et
surfaces minimales}

\begin{lem}[\cite{[De]}]\label{commut}
Soient $f$ et $g$ deux transformations birationnelles sur une
surface $S$ virtuel\-lement isotopes \`a l'identit\'e. Supposons que
$f$ et $g$ commutent~; alors $f$ et $g$ sont simultan\'ement
virtuellement isotopes \`a l'identit\'e.
\end{lem}

\begin{rem}\label{iso}
Un automorphisme $f$ d'une surface $S$
isotope \`a l'identit\'e fixe chaque courbe d'auto-intersection
n\'egative~; pour toute suite de contractions $\psi$ de $S$ vers
un mod\`ele minimal $X$ de $S,$ l'\'el\'ement $\psi
f\psi^{-1}$ est donc un automorphisme de $X$ isotope \`a
l'identit\'e.
\end{rem}

Rappelons que les surfaces minimales sont $\P^1(\C)\times
\P^1(\C)$, $\P^2(\C)$ et les surfaces de \textsc{Hirzebruch}
$\mathsf{F}_n$, $n\geq 2$. Dans des cartes affines bien
choisies $\Aut(\P^1(\C)\times\P^1(\C))$ co\"incide avec
$$(\PGL_2(\C)\times \PGL_2(\C))\rtimes(y,x)$$ et $\Aut(
\mathsf{F}_n)$ se d\'ecrit, pour $n\geq 2$, comme suit 
(\cite{[AS]} chapitre $5$, \cite{[Na]})~:
$$\left\{\left(\frac{\alpha x+P(y)}{(cy+d)^n},
\frac{ay+b}{cy+d}\right)\hs|\hs \left(\begin{array}{cc}
a & b\\
c & d
\end{array}\right)
\in\PGL_2(\C),\hs\alpha\in\C^*,\hs P\in\C[y],\hs\deg P\leq
n\right\}.$$

\section{Sous-groupes nilpotents de groupes d'automorphismes
de surfaces minimales}

Dans la suite du texte, chaque fois que nous dirons que
$\mathsf{N}$ est contenu dans le groupe de \textsc{Jonqui\`eres}
il sera sous-entendu que c'est \`a conjugaison pr\`es.

Soit $\mathsf{N}$ un groupe nilpotent~; posons
$\mathsf{N}^{(0)}:=\mathsf{N}$ et $\mathsf{N}^{(i)}:=
[\mathsf{N},\mathsf{N}^{(i-1)}]$ pour $i\geq 1$.

\begin{pro}\label{surfmin}
Soient $\mathsf{N}$ un sous-groupe n. f. de longueur $k$ de
$\Bir(\P^2(\C))$ et $S$ une surface minimale. Supposons que
$\mathsf{N}^{(k-1)}$ ne soit pas de torsion et soit, \`a
conjugaison birationnelle pr\`es, un sous-groupe de $\Aut(S)$.
Alors $\mathsf{N}$ est, \`a conjugaison birationnelle et indice
fini pr\`es, contenu dans le groupe de \textsc{Jonqui\`eres}.
\end{pro}

Avant de d\'emontrer cette proposition rappelons l'\'enonc\'e
suivant que nous utiliserons \`a plusieurs reprises~:

\begin{thm}[Th\'eor\`eme des r\'epliques de Chevalley]\label{chevalley}
Soient $\mathsf{G}$ un sous-groupe de \textsc{Lie} alg\'ebrique de
$\mathsf{GL}_n(\C)$ et $g$ un \'el\'ement de $G$. Les parties
semi-simple et nilpotente de $g$ appartiennent \`a $\mathsf{G}$.
\end{thm}

\begin{proof}[D\'emonstration]
Remarquons que, puisque $\mathsf{N}$ n'est pas virtuellement
de longueur $k-1$, le groupe $\mathsf{N}^{(k-1)}$ n'est pas fini.
\vs

\noindent {\it 1.} Commen\c{c}ons par consid\'erer le cas o\`u
$S=\P^2(\C)$. Le groupe $\mathsf{N}^{(k-1)}$ est, \`a conjugaison
pr\`es, contenu dans l'un des groupes suivants~:
\begin{align*}
& {\it a. }\hs \{(x+\alpha,y+\beta)\hs|\hs\alpha,\hs\beta \in\C\}
; && {\it c. }\hs \{(\alpha x,y+\beta)\hs|\hs\alpha\in\C^*, \hs
\beta\in \C\}~; \\
&{\it b. }\hs \{(\alpha x,\beta y)\hs|\hs\alpha,\hs\beta \in\C^*\}
;&& {\it d. }\hs \{(\alpha x+\beta y,\alpha y)\hs|\hs\alpha
\in\C^*, \beta\in\C\}.
\end{align*}

\noindent Nous allons consid\'erer ces \'eventualit\'es au cas par
cas.

\vs

\noindent {\it 1.a. } Soient $g=(x+\alpha,y+\beta)$ un \'el\'ement
non trivial de $\mathsf{N}^{(k-1)}$ et $f$ une transformation
birationnelle qui commute \`a $g$. Quitte \`a conjuguer $g$ par
$(y,x)$ nous pouvons supposer que $\alpha\not=0$~; alors $f$
commute, \`a conjugaison pr\`es, \`a $(x+1,y+\beta)$ et est donc
du type $(x+b(y),\nu(y))$. Il en r\'esulte que $\mathsf{N}$ est un
sous-groupe de~: $$\{(x+b(y),\nu(y)) \hs| \hs
b\in\C(y),\hs\nu\in\PGL_2(\C)\}\subset\mathsf{J}.$$

\vs

\noindent {\it 1.b.} Supposons que $\mathsf{N}^{(k-1)}$ contienne
un \'el\'ement $(\alpha x,\beta y)$ non p\'eriodique~; soit $h$
une fonction rationnelle satisfaisant $h(\alpha x,\beta
y)=h(x,y)$. Si $\alpha$ et $\beta$  sont \og non r\'esonants\fg ,
alors $h$ est constante. S'il existe deux entiers $p$ et $q$
premiers entre eux tels que $\alpha^p \beta^q=1,$ la fonction $h$
est invariante par $\overline{\langle(\alpha x,\beta y)
\rangle}^{\hs\mathsf{Z}},$ o\`u l'adh\'erence est prise au sens de
\textsc{Zariski}~; ce groupe est engendr\'e par le flot du champ~:
$$\chi:=qx\frac{\partial} {\partial x}-py \frac{\partial}
{\partial y}.$$ Les niveaux de $h$ sont les trajectoires de $\chi$
donc $h=h(x^py^q)$. Nous en d\'eduisons ce qui suit.

\noindent Soit
$f=(f_1,f_2)$ une transformation birationnelle qui commute \`a
$(\alpha x,\beta y)$~; alors $f$ commute \`a $\overline{
\langle(\alpha x,\beta y)\rangle}^{ \hs\mathsf{Z}}$. 
Si $\alpha$ et $\beta$ sont \og non r\'esonants\fg , $f$ est de la
forme $(\gamma x,\delta y)$ et $\mathsf{N}$ est contenu dans le
groupe de \textsc{Jonqui\`eres}. S'il existe $p$ et $q$ deux
entiers premiers entre eux tels que $\alpha^p\beta^q=1$, alors $f$
commute \`a tous les $(\lambda x,\mu y)$ satisfaisant $\lambda^p
\mu^q=1$, {\it i.e.} $f$ est du type $(xa(x^py^q),yb(x^py^q))$.
Les entiers $p$ et $q$ \'etant premiers entre eux, il existe deux
entiers $u$ et $v$ tels que $uq-pv=1$. \`A conjugaison pr\`es par
$(x^uy^v,x^py^q)$ la transformation $f$ est de la forme
$(c(y)x,\nu(y))~;$ par suite $\mathsf{N}$ est un sous-groupe de
$\mathsf{J}$.

\noindent Supposons que tous les \'el\'ements de $\mathsf{
N}^{(k-1)}$ soient p\'eriodiques~; comme $\mathsf{N}^{
(k-1)}$ n'est pas fini, nous parvenons au m\^eme r\'esultat
en consid\'erant l'adh\'erence de \textsc{Zariski} de
$\mathsf{N}^{(k-1)}$.

\vs

\noindent {\it 1.c.} Le groupe $\mathsf{N}^{(k-1)}$ n'\'etant
pas fini, nous avons l'alternative suivante~:
\begin{itemize}
\item $\mathsf{N}^{(k-1)}$ contient un \'el\'ement de la forme
$(\alpha x,y+\beta)$ avec $\beta\not=0$~;

\item $\{\alpha\in\C^*\hs|\hs(\alpha x,y)\in\mathsf{N}^{(k-1)}\}$
n'est pas fini.
\end{itemize}

Consid\'erons le premier cas. Soit $f$ dans $\mathsf{N}~;$ le
groupe $G=\{g\in\mathsf{N}^{(k-1)}\hs|\hs fg=gf\}$, qui contient
$(\alpha x, y+\beta)$, s'identifie \`a un sous-groupe alg\'ebrique
d'un groupe lin\'eaire. Le th\'eor\`eme des r\'epliques de
\textsc{Chevalley} assure alors que $f$ commute \`a $(x,y
+~\frac{\beta}{\alpha})$~; ainsi $f$ est du type $(\nu(x),y+
b(x))$ et $$\mathsf{N}\subset\{(\nu(x),y+ b(x)) \hs|\hs
\nu\in\PGL_2(\C),\hs b\in\C(x)\}~;$$ quitte \`a conjuguer par
$(y,x)$ nous constatons que $\mathsf{N}$ est inclus dans
$\mathsf{J}$.

Pour la deuxi\`eme \'eventualit\'e, voir {\it 1.b.}

\vs

\noindent {\it 1.d.} Puisque $\mathsf{N}^{(k-1)}$ n'est pas
fini, nous avons l'alternative suivante~:
\begin{itemize}
\item ou bien $\mathsf{N}^{(k-1)}$ contient un \'el\'ement du type
$(\alpha x+\beta y,\alpha y)$ avec $\beta\not=0$~;

\item ou bien $\{\alpha\in\C^*\hs|\hs(\alpha x,\alpha
y)\in\mathsf{N}^{ (k-1)}\}$ n'est pas fini.
\end{itemize}

Consid\'erons la premi\`ere \'eventualit\'e~; le th\'eor\`eme des
r\'epliques de \textsc{Chevalley} assure encore qu'un \'el\'ement
qui commute \`a $(\alpha x+ \beta y,\alpha y)$ commute \`a
$(x+\frac{\beta}{\alpha}y,y)$. Nous en d\'eduisons que
$\mathsf{N}$ est un sous-groupe de $\mathsf{J}$.

Dans le second cas nous obtenons, par {\it 1.b.}, que $\mathsf{N}$
est contenu dans le groupe de \textsc{Jonqui\`eres}.

\vs

\noindent {\it 2.} Etudions le cas o\`u $S=\P^1(\C)
\times\P^1(\C)$. \`A indice fini pr\`es, nous pouvons supposer que
$\mathsf{N}^{(k-1)}$ est un sous-groupe ab\'elien de $\PGL_2(\C)
\times\PGL_2(\C)~;$ il en r\'esulte que $\mathsf{N}^{(k-1)}$ est
inclus, \`a conjugaison pr\`es, dans l'un des groupes suivants~:
$$\{(x+\alpha,y+\beta)\hs|\hs\alpha, \hs\beta \in\C\},
\hspace{6mm} \{(\alpha x,y+\beta)\hs| \hs\alpha\in\C^*, \hs
\beta\in \C\},\hspace{6mm}\{(\alpha x,\beta y)\hs|\hs\alpha,
\hs\beta \in\C^*\}~;$$ le {\it 1.} permet donc de conclure.

\vs

\noindent {\it 3.} Pour finir nous nous int\'eressons au cas
$S=\mathsf{F}_n$, $n\geq 2$. Notons $\pi$ la projection de
$\Aut(\mathsf{F}_n)$ sur $\PGL_2(\C)$~; puisque le groupe $\pi(
\mathsf{N}^{(k-1)})$ est ab\'elien, nous pouvons supposer,
\`a indice fini pr\`es, que $\mathsf{N}^{(k-1)}$ est contenu
dans l'un des groupes suivants~:
\begin{align*}
& {\it a. }\hs \{(\alpha x+P(y),y)\hs|\hs\alpha\in\C^*,\hs P\in
\C[y]\}; &&  \\
& {\it b. }\hs \{(\alpha x+P(y),y+\beta)\hs|\hs\alpha
\in\C^*,\hs\beta\in\C,\hs P\in\C[y]\}; &&  \\
& {\it c. }\hs \{(\alpha x+P(y),\beta y)\hs|\hs\alpha,\hs\beta
\in\C^*,\hs P\in\C[y]\}; &&
\end{align*}
Consid\'erons chacune de ces \'eventualit\'es.\vs

\noindent{\it 3.a.} Supposons que $\mathsf{N}^{(k-1)}$ soit
un sous-groupe de $\{(x+P(y),y)\hs|\hs
P\in\C[y]\}$. Soit $g$ dans $\mathsf{N}^{(k-1)}\setminus
\{\mathsf{id}\}$~; \`a
conjugaison pr\`es, $g$ s'\'ecrit $(x+1,y)$ et {\it 1.a.}
permet de conclure.

\noindent Supposons que $\mathsf{N}^{(k-1)}$ contienne
un \'el\'ement $g$ de la forme $(\alpha x+P(y),y)$ avec $\alpha
\not=1$. Si $\alpha$ est d'ordre fini $q$, alors $g^q$ est
du type $(x+Q(y),y)$ et nous concluons comme pr\'ec\'edemment.
Si $\alpha$ est d'ordre infini, $g$ est conjugu\'e 
\`a $(\alpha x,y)$ et {\it 1.c.} assure que $\mathsf{N}$ est 
inclus dans le groupe de \textsc{Jonqui\`eres}.

 \vs

\noindent{\it 3.b.} Un \'el\'ement $(\alpha x+P(y),y+\beta)$
de $\mathsf{N}^{(k-1)}$, avec $\beta\not=0$, est conjugu\'e \`a
$(\alpha x,y+\beta)$~; par {\it 1.c.} nous obtenons donc
que $\mathsf{N}$ est un sous-groupe de $\mathsf{J}$.\vs

\noindent{\it 3.c.} Si $\{\beta\in\C^*\hs|\hs\beta y\in\pi(
\mathsf{N}^{(k-1)})\}$ est fini, il existe dans $\mathsf{N}^{
(k-1)}$ un \'el\'ement de la forme $(\alpha x+P(y),y)$ et {\it
3.a.} permet de conclure. Supposons donc que $g$ d\'esigne un
\'el\'ement de $\mathsf{N}^{(k-1)}$ du type $(\alpha x+P(y),\beta
y)$ avec $\beta$ d'ordre infini. Si $\beta^i-\alpha\not=0$ pour
tout $0\leq i\leq\deg P$, la transformation $g$ est conjugu\'ee
\`a $(\alpha x,\beta y)$ et, par {\it 1.b.}, nous avons le
r\'esultat annonc\'e. Reste \`a consid\'erer le cas o\`u il existe
un unique $i$ tel que $\beta^i-\alpha=0$ (l'unicit\'e r\'esulte du
fait que $\beta$ est d'ordre infini)~; $g$ est alors conjugu\'e
\`a $(\beta^ix+\gamma y^i,\beta y)$. Il s'en suit qu'une
transformation birationnelle qui commute \`a $g$ commute \`a
$(\beta^ix,\beta y)$ (th\'eor\`eme des r\'epliques de
\textsc{Chevalley} \og adapt\'e\fg )~;
nous concluons \`a l'aide de {\it 1.b.}
\end{proof}

\section{Dynamique et groupes nilpotents de g\'en\'eration
finie}

\begin{lem}\label{degdyn1}
Soient $\mathsf{N}$ un groupe nilpotent non ab\'elien de longueur
de nilpotence $k$ et $\rho$ un morphisme injectif de $\mathsf{N}$
dans $\Bir(\P^2(\C)).$ Nous sommes dans l'une des situations
suivantes~:
\begin{itemize}
\item ou bien pour tout $g$ dans $\rho(\mathsf{N})$, nous avons
$\lambda(g)=1$~;

\item ou bien $\mathsf{N}^{(k-1)}$ est de torsion.
\end{itemize}
\end{lem}

\begin{proof}[D\'emonstration]
Soient $\mathsf{f}$ dans $\mathsf{N}$ et $\mathsf{g}$ dans
$\mathsf{N}^{(k-2)}~;$ consid\'erons l'\'el\'ement $\mathsf{h}$
d\'efini par $\mathsf{h}=[\mathsf{f},\mathsf{g}].$ Supposons que
$\mathsf{h}$ soit non p\'eriodique. Puisque $\mathsf{N}$ est de
longueur de nilpotence $k$, nous avons $[\mathsf{f},\mathsf{h}]=
[\mathsf{g},\mathsf{h}]=\mathsf{id}.$ La remarque \ref{dist} et
le lemme \ref{degdyn} assurent que le premier degr\'e dynamique de
$\rho(\mathsf{h})$ vaut $1.$ Montrons que pour tout $\ell$ dans
$\mathsf{N}$, nous avons $\lambda(\rho(\ell))=1$. La
transformation $\rho(\mathsf{h})$ satisfait l'alternative suivante
(\cite{[Di-Fa]})~:
\begin{itemize}
\item $\rho(\mathsf{h})$ pr\'eserve une unique fibration
rationnelle ou elliptique~;

\item $\rho(\mathsf{h})$ est virtuellement isotope \`a
l'identit\'e.
\end{itemize}

Dans le premier cas, $\rho(\ell)$ laisse, par commutation, une
fibration invariante donc $\rho(\ell)$ est de premier
degr\'e dynamique $1$. Consid\'erons la seconde \'eventualit\'e~;
pla\c{c}ons nous sur la surface $S$ o\`u $\rho(\mathsf{h})$ est un
automorphisme et notons encore $\rho(\mathsf{h})$ (resp. $\rho(
\ell)$) le conjugu\'e de $\rho(\mathsf{h})$ (resp. $\rho( \ell)$).
Supposons $\lambda(\rho(\ell))>1$. Notons $G=\overline{\langle
\rho(\mathsf{h})\rangle}^{\hs\mathsf{Z}}\subset\Aut(S)~;$ comme
$\rho(\mathsf{h})$ est d'ordre infini $G$ est un groupe de
\textsc{Lie} ab\'elien qui commute \`a $\rho(\ell)$~: ceci
contredit $\lambda(\rho(\ell))>1.$\footnote{\hs Ce dernier
argument est d\^u \`a S. \textsc{CANTAT}.}

Supposons que tous les commutateurs du type $[f,g]$ o\`u
$f$ d\'esigne un \'el\'ement de $\rho(\mathsf{N})$ et $g$
un \'el\'ement de $\rho(\mathsf{N}^{(k-2)})$ soient
p\'eriodiques~; alors $\rho(\mathsf{N}^{(k-1)})$ \'etant
ab\'elien, il est de torsion.
\end{proof}

\begin{egs}
{\it 1. } Posons $H_0=\langle (-x,-y) \rangle$ et pour
tout $j\geq 1$~: $$\alpha_j=\beta_j=\exp(\mathrm{i}
\pi/2^j),\hspace{6mm} g_j=(y,c_jy^{2j+1}+x),$$
$$\varphi_j=g_1^2\circ \ldots\circ g_j^2,\hspace{6mm}
H_j=\varphi_j\langle (\alpha_j x,\beta_j y)\rangle\varphi_j^{-1}$$
o\`u $c_j$ appartient \`a $\C^*.$ Alors $H=\cup_{j\geq 0} H_j$ est
un groupe ab\'elien d\'enombrable dont tous les \'el\'ements sont p\'eriodiques
(\cite{[La], [Wr]}).\vs

{\it 2.} Le groupe $\mathsf{N}=\langle(\exp(-2\mathrm{i}\pi/p)x,y)
,(x,xy),(x,\exp(2\mathrm{i}\pi/p)y)\rangle$ est nilpotent de longueur
$2$ et $\mathsf{N}^{(1)}=\langle(x,\exp(2\mathrm{i}\pi/p)y)
\rangle$ est de torsion.
\end{egs}

\begin{lem}\label{fibvirtiso}
Soient $\mathsf{N}$ un groupe n. f. de longueur $k>1$ et $\rho$ 
un morphisme injectif de $\mathsf{N}$ dans
$\Bir(\P^2(\C)).$ Supposons que $\mathsf{N}^{(k-1)}$ ne soit pas
de torsion. Le groupe $\rho(\mathsf{N})$ satisfait l'une des
propri\'et\'es suivantes~:
\begin{itemize}
\item ou bien $\rho(\mathsf{N})$ pr\'eserve une fibration 
rationnelle ou elliptique~;

\item ou bien les \'el\'ements de $\rho(\mathsf{N}^{(k-1)})$ sont
virtuellement isotopes \`a l'identit\'e.
\end{itemize}
\end{lem}

\begin{proof}[D\'emonstration]
Soit $h$ un \'el\'ement non trivial de $\mathsf{N}^{(k-1)}$~;
le lemme \ref{degdyn1} assure que
$\lambda(\rho(h))=1$. D'apr\`es \cite{[Di-Fa]} ou bien
$\rho(h)$ est virtuellement isotope \`a l'identit\'e, ou bien
$\rho(h)$ pr\'eserve une unique fibration $\mathcal{F}$
rationnelle ou elliptique. Pla\c{c}ons nous dans cette 
seconde \'eventualit\'e~; comme $\rho(\mathsf{N})$ commute 
\`a $\rho(h)$, tout \'el\'ement de $\rho(\mathsf{N})$ 
pr\'eserve $\mathcal{F}$ (c'est l'unicit\'e).
\end{proof}

Rappelons que tout sous-groupe d'un
groupe nilpotent de g\'en\'eration finie est nilpotent de
g\'en\'eration finie (\cite{[Ma]}, th\'eor\`eme 9.16).

\begin{pro}\label{blablanil2}
Soit $\mathsf{N}$ un groupe n. f. de longueur $k$ non ab\'elien et
de g\'en\'eration finie. Soit $\rho$ un morphisme injectif de
$\mathsf{N}$ dans $\Bir(\P^2(\C)).$ Supposons que
$\rho(\mathsf{N}^{(k-1)})$ ne soit pas de torsion et que ses
\'el\'ements soient virtuellement isotopes \`a l'identi\-t\'e~;
alors $\rho(\mathsf{N})$ est, \`a indice fini et conjugaison
pr\`es, contenu dans le groupe de \textsc{Jonqui\`eres}. 
\end{pro}

\begin{proof}[D\'emonstration]
Le groupe $\rho(\mathsf{N}^{(k-1)})$ est de g\'en\'eration
finie~; notons $\{a_1,\ldots,a_n\}$ une partie g\'en\'eratrice de
$\rho(\mathsf{N}^{(k-1)})$. Comme $\rho(\mathsf{N}^{(k-1)})$ est
ab\'elien le lemme \ref{commut} assure que les $a_i$ sont
simultan\'ement virtuellement isotopes \`a l'identit\'e. D'apr\`es
la remarque \ref{iso} le groupe $\langle a_1^q,\ldots,
a_n^q\rangle$ est, \`a conjugaison pr\`es, contenu dans le
groupe d'automorphismes d'une surface minimale~; la proposition
\ref{surfmin} permet de conclure.
\end{proof}

\begin{pro}\label{torsion}
Soient $\mathsf{N}$ un groupe n. f. de longueur $k>1$ 
de g\'en\'eration finie et
$\rho$ un morphisme injectif de $\mathsf{N}$ dans $\Bir(\P^2(
\C))$. Supposons que $\rho(\mathsf{N}^{(k-1)})$ soit de torsion.
Alors $\rho(\mathsf{N})$ v\'erifie l'une des propri\'et\'es
suivantes~:
\begin{itemize}
\item $\rho(\mathsf{N})$ laisse une fibration elliptique ou
rationnelle invariante~;

\item $\rho(\mathsf{N})$ est fini.
\end{itemize}
\end{pro}

\begin{proof}[D\'emonstration]
Le groupe $\mathsf{N}^{ (k-1)}$ est un groupe ab\'elien, de
g\'en\'eration finie, dont tous les \'el\'ements sont
p\'eriodiques donc fini. Il existe une surface $M$ et une
transformation birationnelle $\eta\hs\colon \hs\P^2(\C)\dasharrow
M$ telles que $\eta\rho( \mathsf{N}^{(k-1)}) \eta^{-1}$ soit
inclus dans $\Aut(M)$. Soit $S$ le quotient de $M$ par $\eta\rho(
\mathsf{N}^{(k-1)})\eta^{-1}$~; notons $\tilde{S}$ la
d\'esingularis\'ee de $S$. La surface $\tilde{S}$ est
unirationnelle donc rationnelle (\emph{voir} \cite{[Be]}). Par
suite $\mathsf{N}_1:=\mathsf{N}/\mathsf{N}^{(k-1)}$ se plonge dans
$\Bir(\P^2(\C))$~; notons $\rho_1$ ce plongement. Le groupe
$\mathsf{N}_1$ est nilpotent de longueur de nilpotence $k-1$.
D'apr\`es le lemme \ref{degdyn1} nous avons l'alternative
suivante~:

\vs

\begin{itemize}
\item[{\it a.}] pour tout $g$ dans $\mathsf{N}_1$ nous avons
$\lambda(\rho_1(g))=1$~;

\item[{\it b.}] $\mathsf{N}_1^{(k-2)}$ est de torsion.
\end{itemize}

\vs

\noindent {\it a.} D'apr\`es le lemme \ref{fibvirtiso} et 
la proposition \ref{blablanil2}, le groupe
$\rho_1(\mathsf{N}_1)$ pr\'eserve une fibration elliptique
ou rationnelle~; montrons qu'il en est de m\^eme
pour $\rho(\mathsf{N})$. Notons
$p$ la projection de $M$ sur $S$ et $F$ la
fibration invariante par $\rho_1(\mathsf{N}_1)$. La fibration $\mathcal{F}=F\circ p$ est invariante
par $\rho(\mathsf{N})$~; \`a tout \'el\'ement $f$ de $\rho(
\mathsf{N})$ nous pouvons donc associer $\nu_f$ dans
$\Aut(\P^1(\C))$ satisfaisant~: $\mathcal{F}\circ
f=\nu_f\circ\mathcal{F}$. Soit $\nu$ le morphisme d\'efini
par~: \begin{eqnarray}
\nu\hs\colon\hs\rho(\mathsf{N})&\to&\PGL_2(\C)\nonumber\\
f&\mapsto&\nu_f\nonumber
\end{eqnarray}
Comme $\rho(\mathsf{N})$ n'est pas virtuellement ab\'elien,
$\ker\nu$ n'est pas fini~; la fibre g\'en\'erique de $\mathcal{F}$ a donc
un groupe d'automorphismes non fini.
Il en r\'esulte que $\mathcal{F}$ est rationnelle ou
elliptique. 

\vs

\noindent {\it b.} Si $\mathsf{N}_1^{(k-2)}$ est de torsion,
$\mathsf{N}_2= \mathsf{N}_1/ \mathsf{N}_1^{(k-2)}$ se plonge dans
$\Bir(\P^2(\C))$~; notons $\rho_2$ ce plongement. Nous pouvons
alors reprendre le m\^eme raisonnement... ainsi~:
\begin{itemize}
\item ou bien il existe $1\leq i\leq k-1$ tel que $\rho_i(
\mathsf{N}_i)$ soit virtuellement contenu dans le groupe de
\textsc{Jonqui\`eres} et, par suite, $\rho(\mathsf{N})$ 
l'est aussi~;

\item ou bien $\mathsf{N}_{k-1}$ est de torsion donc
fini et $\mathsf{N}$ aussi.
\end{itemize}
\end{proof}

\section{Conclusion}

\begin{rem}\label{j0}
Un sous-groupe nilpotent non fini de $\PGL_2(\mathbf{k})$, avec
$\mathbf{k}=\C$ ou $\C(y)$, est n\'ecessairement virtuellement
ab\'elien.
\end{rem}

Avant de donner une preuve du th\'eor\`eme \ref{nilpotent}
nous \'etablissons un \'enonc\'e pour les
sous-groupes nilpotents, non ab\'eliens et de g\'en\'eration
finie de $\Bir(\P^2(\C))$.

\begin{thm}\label{nilgenfini}
Soient $\mathsf{N}$ un groupe n. f. de longueur $k$ non
ab\'elien, de g\'en\'eration fi\-nie et $\rho$ un morphisme
injectif de $\mathsf{N}$ dans $\Bir(\P^2(\C))$. Le groupe
$\rho(\mathsf{N})$ v\'erifie l'une des propri\'et\'es suivantes~:
\begin{itemize}
\item $\rho(\mathsf{N})$ pr\'eserve une fibration elliptique
ou rationnelle~;

\item $\rho(\mathsf{N})$ est fini.
\end{itemize}
\end{thm}

\begin{proof}[D\'emonstration]
L'\'enonc\'e d\'ecoule des lemmes 
\ref{degdyn1}, \ref{fibvirtiso}, des propositions 
\ref{blablanil2} et \ref{torsion}. 
\end{proof}

Nous allons, pour finir, d\'emontrer le th\'eor\`eme
\ref{nilpotent}.

\begin{proof}[D\'emonstration du th\'eor\`eme \ref{nilpotent}]
Notons $k$ la longueur de nilpotence de $\rho(\mathsf{N})$.
Consi\-d\'erons $\Sigma$ l'ensemble des sous-groupes n. f. de
longueur $k$ et de g\'en\'eration finie de $\rho(\mathsf{N})$.
Si tous les \'el\'ements de $\Sigma$ sont finis alors 
$\rho(\mathsf{N})$ est de torsion~; sinon $\Sigma$
contient un \'el\'ement $\mathsf{G}$ qui ne soit pas fini.
D'apr\`es le th\'eor\`eme \ref{nilgenfini}, $\rho(\mathsf{G})$ 
pr\'eserve une fibration elliptique ou rationnelle $\mathcal{F}.$
Dans les deux cas les \'el\'ements de $\mathsf{G}^{(k
-1)}$ pr\'eservent $\mathcal{F}$ fibre \`a fibre. Soit
$f$ dans $\mathsf{G}^{(k-1)}$~; en utilisant la commutation de $f$ \`a
$\rho(\mathsf{N})$ nous observons l'alternative suivante~:
\begin{itemize}
\item $\rho(\mathsf{N})$ laisse aussi la fibration $\mathcal{F}$
invariante,

\item $f$ pr\'eserve deux fibrations distinctes fibre \`a fibre.
\end{itemize}

On constate que dans le premier cas $\mathsf{N}$ est 
virtuellement m\'etab\'elien. En effet un groupe de 
transformations birationnelles de $\P^2(\C)$ qui pr\'eservent
une fibration elliptique est virtuellement m\'etab\'elien.
Si $\mathcal{F}$ est rationnelle, alors $\rho(\mathsf{N})$
est, \`a conjugaison pr\`es, un sous-groupe de $\mathsf{J}$. 
Notons $\pi$ la projection de $\mathsf{J}$ sur $\PGL_2(\C)~;$ le
groupe $\pi(\rho(\mathsf{N}))$ est un sous-groupe nilpotent
de $\PGL_2(\C)$, il est donc virtuellement ab\'elien 
(remarque \ref{j0})~: nous pouvons donc supposer que 
$\pi(\rho(\mathsf{N}^{(i)}))=\{\mathsf{id}\}$ pour $1\leq
i\leq k.$ En particulier $\mathsf{N}^{(1)}$ est un sous-groupe
nilpotent de $\PGL_2(\C(y))$, il est donc virtuellement 
ab\'elien (toujours par la remarque \ref{j0}). 

Dans la seconde \'eventualit\'e $f$ est p\'eriodique. Si 
cela a lieu pour tout choix de $f$ dans $\mathsf{G}^{(k-1)}$ 
alors $\mathsf{G}$ est virtuellement de longueur $k-1$.
\end{proof}

\begin{eg}
Les groupes suivants sont des sous-groupes nilpotents non ab\'eliens
et non finis de $\mathsf{J}$~:
$$\langle(x+\alpha\beta,y),\hs(x+\alpha y,y),\hs(x,y+\beta)\rangle ,$$
$$\langle(x+1,y),\hs(x+y,y),\hs(x+a(y),y-1)\rangle$$
o\`u $\alpha$, $\beta$ appartiennent \`a $\C^*$ et
$a$ \`a $\C(y)$.
\end{eg}

\end{document}